\documentclass[12pt]{article}
\usepackage[utf8]{inputenc}
\usepackage[T1]{fontenc}

\usepackage[]{amsfonts}
\usepackage{amssymb}
\usepackage{amsmath}
\def\couleur(#1 #2 #3)
	{
\special{ps::[end]}}

\def\bx#1{\setbox1=\hbox{\kern3pt{#1}\kern3pt}			
 \dimen1=\ht1 \advance\dimen1 by 3pt \dimen2=\dp1 \advance\dimen2 by 3pt
 \setbox1=\hbox{\vrule height\dimen1 depth\dimen2\box1\vrule}%
 \setbox1=\vbox{\hrule\box1\hrule}%
 \advance\dimen1 by .4pt \ht1=\dimen1
 \advance\dimen2 by .4pt \dp1=\dimen2 \box1\relax}

\def\wbb#1{\kern#1em}
\def\vci{\vrule  width.02em height1.47ex depth-.0ex}		
\def\11{{\rm\wbb{.2}\vci\wbb{-.37}1}}

\def\underset#1#2{\mathrel{\mathop{\kern0pt #2}\limits_{#1}}}

\def\overset#1#2{\mathrel{\mathop{\kern0pt #2}\limits^{#1}}}

\def\Supp{\mathop{\rm Supp}\nolimits}

\newtheorem{Thrm}{Theorem}[section]

\begin{document}

\title{H\"ormander's solution of the  $\bar \partial $ -equation with compact support.}

\author{Eric Amar}

\date{ }
\maketitle
 \ \par 
\ \par 
\ \par 
\renewcommand{\abstractname}{Abstract}

\begin{abstract}
This work is a complement of the study on H\"ormander's solution
 of the  $\displaystyle \bar \partial $  equation initialised
 by H. Hedenmalm. Let  $\varphi $  be a strictly plurisubharmonic
 function of class  $\displaystyle {\mathcal{C}}^{2}$  in  $\displaystyle
 {\mathbb{C}}^{n}$  ; let  $\displaystyle c_{\varphi }(z)$  be
 the smallest eigenvalue of  $\displaystyle i\partial \bar \partial
 \varphi $  then  $\displaystyle \forall z\in {\mathbb{C}}^{n},\
 c_{\varphi }(z)>0.$   We denote by  $\displaystyle L_{p,q}^{2}({\mathbb{C}}^{n},e^{\varphi
 })$   the  $\displaystyle (p,q)$  currents with coefficients
 in  $\displaystyle L^{2}({\mathbb{C}}^{n},e^{\varphi }).$ \ \par 
\quad We prove that if  $\displaystyle \omega \in L^{2}_{p,q}({\mathbb{C}}^{n},e^{\varphi
 }),\ \bar \partial \omega =0$  for  $\displaystyle q<n$  then
 there is a solution  $\displaystyle u\in L^{2}_{p,q-1}({\mathbb{C}}^{n},c_{\varphi
 }e^{\varphi })$  of  $\displaystyle \bar \partial u=\omega .$ \ \par 
\quad This is done via a theorem giving a solution with compact support
 if the data has compact support.\ \par 
\end{abstract}

\section{Introduction.}
\quad Let  $\varphi $  be a strictly plurisubharmonic function of class
  ${\mathcal{C}}^{2}$  in a domain  $\displaystyle \Omega $ 
 in  $\displaystyle {\mathbb{C}}^{n}.$  Let  $\displaystyle c_{\varphi
 }(z)$  be the smallest eigenvalue of  $i\partial \bar \partial
 \varphi (z),$  then  $\displaystyle \forall z\in \Omega ,\ c_{\varphi
 }(z)>0.$   We denote by  $\displaystyle L^{2,c}_{p,q}(\Omega
 ,e^{\varphi })$  the currents in  $\displaystyle L^{2}_{p,q}(\Omega
 ,e^{\varphi })$  with compact support in  $\displaystyle \Omega
 $  and   $\displaystyle {\mathcal{H}}_{p}(\Omega ,e^{-\varphi
 })$  the space of all  $(p,\ 0),\ \bar \partial $  closed forms
 in  $\displaystyle L^{2}(\Omega ,\ e^{-\varphi }).$ \ \par 
\ \par 
\quad In the case  $\displaystyle n=1,$  H. Hedenmalm~\cite{Hedenmalm13}
 proved the following theorem\ \par 

\begin{Thrm}
Let  $\displaystyle \omega \in L^{2}_{0,1}({\mathbb{C}},\ e^{\varphi
 })$  then there is  $\displaystyle u\in L^{2}({\mathbb{C}},c_{\varphi
 }e^{\varphi })$  such that  $\displaystyle \bar \partial u=\omega ,$  and\par 
\quad \quad \quad \quad \quad \quad \quad  $\displaystyle \ {\left\Vert{u}\right\Vert}_{L^{2}({\mathbb{C}},c_{\varphi
 }e^{\varphi })}\leq C{\left\Vert{\omega }\right\Vert}_{L^{2}({\mathbb{C}},e^{\varphi
 })},$ \par 
provided that  $\displaystyle \omega \perp {\mathcal{H}}_{0}({\mathbb{C}},\
 e^{-\varphi }).$  
\end{Thrm}
\quad As he prognosticated  in his remark 1.3 in~\cite{Hedenmalm13}
 this theorem was generalised to Stein manifolds in~\cite{AmarHeden13},
 and for domains in  $\displaystyle {\mathbb{C}}^{n}$  we got :\ \par 

\begin{Thrm}
~\label{cS6}Let  $\displaystyle \Omega $  be a pseudo convex
 domain in  ${\mathbb{C}}^{n}$ ; if  $\omega \in L^{2,c}_{p,q}(\Omega
 ,\ e^{\varphi })$  with  $\displaystyle \bar \partial \omega
 =0$  if  $\displaystyle q<n$  and  $\displaystyle \omega \in
 L^{2}_{p,q}(\Omega ,e^{\varphi })$  with  $\displaystyle \omega
 \perp {\mathcal{H}}_{n-p}(\Omega ,\ e^{-\varphi })$  if  $\displaystyle
 q=n,$  then there is  $\displaystyle u\in L^{2}_{p,q-1}(\Omega
 ,c_{\varphi }e^{\varphi })$  such that  $\displaystyle \bar
 \partial u=\omega ,$  and		      	 $\displaystyle \ {\left\Vert{u}\right\Vert}_{L^{2}(\Omega
 ,c_{\varphi }e^{\varphi })}\leq C{\left\Vert{\omega }\right\Vert}_{L^{2}(\Omega
 ,e^{\varphi })}.$  
\end{Thrm}
\quad Still pushed by H. Hedenmalm in several e-mails, the aim  is
 to improve this theorem in the special case of  $\displaystyle
 {\mathbb{C}}^{n}$  by removing the assumption on the compact
 support for  $\omega .$  Precisely we shall prove the following\ \par 

\begin{Thrm}
Let  $\displaystyle \omega \in L^{2}_{p,q}({\mathbb{C}}^{n},\
 e^{\varphi })$  with  $\displaystyle \bar \partial \omega =0$
  if  $\displaystyle q<n$  and with  $\displaystyle \omega \perp
 {\mathcal{H}}_{n-p}(\Omega ,\ e^{-\varphi })$  if  $\displaystyle
 q=n,$  then there is  $\displaystyle u\in L^{2}_{p,q-1}({\mathbb{C}}^{n},c_{\varphi
 }e^{\varphi })$  such that  $\displaystyle \bar \partial u=\omega ,$  and\par 
\quad \quad \quad \quad \quad \quad \quad  $\displaystyle \ {\left\Vert{u}\right\Vert}_{L^{2}({\mathbb{C}}^{n},c_{\varphi
 }e^{\varphi })}\leq C{\left\Vert{\omega }\right\Vert}_{L^{2}({\mathbb{C}}^{n},e^{\varphi
 })}.$ 
\end{Thrm}
I don't know if this is still true for a proper pseudo convex
 domain  $\displaystyle \Omega $  of  $\displaystyle {\mathbb{C}}^{n}$
  ; the method I have works only for  $\displaystyle {\mathbb{C}}^{n}.$ \ \par 
\quad The proof is strongly based on theorem~\ref{cS6} via the theorem~\ref{cS1},
 which can be seen as an auto improvement of theorem~\ref{cS6}.\ \par 

\section{Solutions with compact support.}
\quad The first aim is to prove the following theorem which is an improvement
 of theorem~\ref{cS6}, because here the solution is also with
 compact support.\ \par 

\begin{Thrm}
~\label{cS1}Let  $\displaystyle \Omega $  be a pseudo convex
 domain in  ${\mathbb{C}}^{n}$ and   $\varphi $  be a strictly
 plurisubharmonic function of class  ${\mathcal{C}}^{2}$  in
  $\displaystyle \Omega .$  Let  $\displaystyle c_{\varphi }(z)$
  be the smallest eigenvalue of  $\partial \bar \partial \varphi (z).$ \par 
There is a constant  $\displaystyle C>0$  such that\par 
\quad $\displaystyle if\ q<n,\ \forall \omega \in L^{2,c}_{p,q}(\Omega
 ,\ e^{\varphi }),\ \bar \partial \omega =0,\ \exists u\in L^{2,c}_{p,q-1}(\Omega
 ,c_{\varphi }e^{\varphi }),\ \bar \partial u=\omega ,\ {\left\Vert{u}\right\Vert}\leq
 C{\left\Vert{\omega }\right\Vert}.$ \par 
If  $\displaystyle q=n,$  we make the hypothesis that  $\displaystyle
 \omega \perp {\mathcal{H}}_{n-p}(\Omega ,\ e^{-\varphi })$ 
 and we still get :\par 
\quad \quad \quad $\displaystyle \forall \omega \in L^{2,c}_{p,n}(\Omega ,\ e^{\varphi
 }),\ \ \exists u\in L^{2,c}_{p,n-1}(\Omega ,c_{\varphi }e^{\varphi
 }),\ \bar \partial u=\omega ,\ {\left\Vert{u}\right\Vert}\leq
 C{\left\Vert{\omega }\right\Vert}.$ 
\end{Thrm}
\quad Proof.\ \par 
The idea is the same we used in~\cite{AnGrauLr12}: we add a weight
 to force the solution to be small outside a fixed compact set,
 then we take a limit.\ \par 
\quad Let  $\displaystyle \omega \in L^{2,c}_{p,q}(\Omega ,\ e^{\varphi
 })$  with  $\displaystyle \bar \partial \omega =0$  if  $\displaystyle
 q<n$  and  $\displaystyle \omega \in L^{2}_{p,q}(\Omega ,e^{\varphi
 })$  with  $\displaystyle \omega \perp {\mathcal{H}}_{n-p}(\Omega
 ,\ e^{-\varphi })$  if  $\displaystyle q=n.$  Take a p.c. domain
  $\displaystyle D$  such that  $\displaystyle \bar D\subset
 \Omega ,$  containing the support of  $\omega $  and defined
 by a p.s.h. smooth function  $\rho ,$  i.e.  $\displaystyle
 D:=\lbrace z\in \Omega ::\rho (z)<0\rbrace .$  This is always
 possible because  $\displaystyle \Omega $  is Stein. Fix  $\displaystyle
 \epsilon >0$  and set  $\displaystyle D_{\epsilon }:=\lbrace
 z\in \Omega ::\rho (z)<\epsilon \rbrace ,$  then  $\displaystyle
 D\subset D_{\epsilon }.$ \ \par 
Let  $\displaystyle \chi \in {\mathcal{C}}^{\infty }({\mathbb{R}})$
  such that  $\chi $  is increasing and convex and such that :\ \par 
\quad \quad \quad \quad $\displaystyle \forall t\leq 0,\ \chi (t)=0\ ;\ \forall t>0,\
 \chi '(t)>0,\ \chi "(t)>0.$ \ \par 
The function  $\sigma :=\chi (\rho )$  is also p.s.h. and smooth.
 We have that  $\displaystyle \sigma >\chi (\epsilon )>0$  outside
  $\bar D_{\epsilon }$  so, with  $\displaystyle \psi _{k}:=\varphi
 +k\sigma $  we have that  $\displaystyle \partial \bar \partial
 \psi _{k}=\partial \bar \partial \varphi +k\partial \bar \partial
 \sigma $  is still positive definite, so  $\psi _{k}$  is s.p.s.h.
 and "big" outside  $\bar D_{\epsilon }.$  Moreover in  $D$ 
 we have  $\displaystyle \psi _{k}=\varphi .$  So we apply theorem~\ref{cS6}
 with  $\psi _{k}$  and we get the existence of a solution  $\displaystyle
 u_{k}\in L^{2}_{p,q-1}(\Omega ,c_{\psi _{k}}e^{\psi _{k}})$
  of  $\displaystyle \bar \partial u_{k}=\omega ,$  with\ \par 
\quad \quad \quad \begin{equation}  \ {\left\Vert{u_{k}}\right\Vert}_{L^{2}(\Omega
 ,c_{\psi _{k}}e^{\psi _{k}})}\leq C{\left\Vert{\omega }\right\Vert}_{L^{2}(\Omega
 ,e^{\varphi })},\label{cS2}\end{equation}\ \par 
because on  $\displaystyle \Supp \ \omega $  we have  $\displaystyle
 \forall k\in {\mathbb{N}},\ \psi _{k}=\varphi .$ \ \par 
But, because  $\sigma $  is  p.s.h., we have that  $\displaystyle
 c_{\psi _{k}}\geq c_{\varphi }$  and  $\displaystyle \forall
 z\notin \bar D,\ e^{\psi _{k}(z)}=e^{k\sigma (z)}e^{\varphi
 (z)}$  with  $\displaystyle \sigma (z)>0,$  this means that
  $\displaystyle u_{k}$  must be small in  $\displaystyle \bar
 D_{\epsilon }^{c}.$  \ \par 
Now we proceed as in~\cite{AnGrauLr12}. We have by~(\ref{cS2}) :\ \par 
\quad \quad \quad $\displaystyle \ \int_{\Omega }{\left\vert{u_{k}}\right\vert
 ^{2}c_{\varphi }e^{\psi _{k}}dm}\leq C{\left\Vert{\omega }\right\Vert}_{L^{2}(\Omega
 ,e^{\varphi })}$ \ \par 
because  $\displaystyle c_{\psi _{k}}\geq c_{\varphi }.$  Then,
 by  $\displaystyle e^{\psi _{k}}=e^{k\sigma }e^{\varphi }\geq
 e^{k\chi (\epsilon )}e^{\varphi }$  on  $\displaystyle D_{\epsilon
 }^{c},$  we get that\ \par 
\quad \quad \quad $\displaystyle e^{k\chi (\epsilon )}\int_{\Omega \backslash D_{\epsilon
 }}{\left\vert{u_{k}}\right\vert ^{2}c_{\varphi }e^{\varphi }dm}\leq
 \int_{\Omega \backslash D_{\epsilon }}{\left\vert{u_{k}}\right\vert
 ^{2}c_{\varphi }e^{\varphi }e^{k\sigma }dm}\leq C{\left\Vert{\omega
 }\right\Vert}_{L^{2}(\Omega ,e^{\varphi })}^{2}$ \ \par 
hence\ \par 
\quad \quad \quad \begin{equation}  \ \int_{\Omega \backslash D_{\epsilon }}{\left\vert{u_{k}}\right\vert
 ^{2}c_{\varphi }e^{\varphi }dm}\leq Ce^{-k\chi (\epsilon )}{\left\Vert{\omega
 }\right\Vert}_{L^{2}(\Omega ,e^{\varphi })}^{2}.\label{cS3}\end{equation}\
 \par 
\quad On the other hand we have\ \par 
\quad \quad \quad $\displaystyle \ \int_{\Omega }{\left\vert{u_{k}}\right\vert
 ^{2}c_{\varphi }e^{\varphi }dm}\leq C{\left\Vert{\omega }\right\Vert}_{L^{2}(\Omega
 ,e^{\varphi })}$ \ \par 
hence the sequence  $\displaystyle \lbrace u_{k}\rbrace _{k\in
 {\mathbb{N}}}$  is uniformly bounded in  $\displaystyle L^{2}_{(p,q-1)}(\Omega
 ,c_{\varphi }e^{\varphi }).$  So there is a subsequence, still
 denoted  $\displaystyle \lbrace u_{k}\rbrace _{k\in {\mathbb{N}}},$
  converging weakly to  $\displaystyle u\in L^{2}_{(p,q-1)}(\Omega
 ,c_{\varphi }e^{\varphi }),$  i.e.\ \par 
\quad \quad \quad $\displaystyle \forall f\in L^{2}_{n-p,n-q+1}(\Omega ,c_{\varphi
 }e^{\varphi }),\ {\left\langle{u_{k},f}\right\rangle}\rightarrow
 {\left\langle{u,f}\right\rangle}.$ \ \par 
\quad To see that this form  $u$  is  $0\ a.e.$  on  $\Omega \backslash
 D_{\epsilon }$  let us take a component  $u_{I,J}$  of it ;
 it is the weak limit of the sequence of functions  $\lbrace
 u_{k,I,J}\rbrace $  which means, with the notations  $v:=u_{I,J},\
 v_{k}:=u_{k,I,J},\ d\mu :=c_{\varphi }e^{\varphi }dm,$ \ \par 
\quad \quad \quad \quad \quad  $\displaystyle \forall g\in L^{2}(\Omega ,c_{\varphi }e^{\varphi
 }),\ \int_{\Omega }{vgd\mu }=\lim _{k\rightarrow \infty }\int_{\Omega
 }{v_{k}gd\mu }.$ \ \par 
As usual take  $\displaystyle g:=\frac{\bar v}{\left\vert{v}\right\vert
 }{\11}_{E}$  where  $E:=\lbrace \left\vert{v}\right\vert >0\rbrace
 \cap (\Omega \backslash D_{\epsilon })$  then we get\ \par 
\quad \quad \quad \quad 	\begin{equation}  \ \int_{\Omega }{vgd\mu }=\int_{E}{\left\vert{v}\right\vert
 d\mu }=\lim _{k\rightarrow \infty }\int_{\Omega }{v_{k}gd\mu
 }=\lim _{k\rightarrow \infty }\int_{E}{\frac{v_{k}\bar v}{\left\vert{v}\right\vert
 }d\mu }.\label{cS4}\end{equation}\ \par 
Now we have by Cauchy-Schwarz\ \par 
\quad \quad \quad \quad \quad  $\displaystyle \ \left\vert{\ \int_{E}{\frac{v_{k}\bar v}{\left\vert{v}\right\vert
 }d\mu }}\right\vert \leq {\left\Vert{v_{k}}\right\Vert}_{L^{2}(E,d\mu
 )}{\left\Vert{{\11}_{E}}\right\Vert}_{L^{2}(E,d\mu )}.$ \ \par 
But\ \par 
\quad \quad \quad \quad \quad  $\displaystyle \ {\left\Vert{v_{k}}\right\Vert}_{L^{2}(E,d\mu
 )}^{2}\leq \int_{\Omega \backslash D_{\epsilon }}{\left\vert{u_{k}}\right\vert
 ^{2}d\mu }\leq Ce^{-k\chi (\epsilon )}{\left\Vert{\omega }\right\Vert}_{L^{2}(\Omega
 ,e^{\varphi })}^{2}$ \ \par 
by~(\ref{cS3}) because  $\displaystyle \chi (\epsilon )>0.$ \ \par 
Hence, by~(\ref{cS4}),\ \par 
\quad \quad \quad \quad \quad  $\displaystyle \ \left\vert{\ \int_{E}{\left\vert{v}\right\vert
 d\mu }}\right\vert ^{2}\leq C{\left\Vert{{\11}_{E}}\right\Vert}_{L^{2}(E,d\mu
 )}\lim \ _{k\rightarrow \infty }\ e^{-k\chi (\epsilon )}{\left\Vert{\omega
 }\right\Vert}_{L^{2}(\Omega ,e^{\varphi })}^{2}=0$ \ \par 
which implies  $\mu (E)=m(E)=0$  because on  $E,\ \left\vert{v}\right\vert
 >0.$ \ \par 
\quad This being true for all components of  $u,$  we get that the
 form  $u$  is  $0\ a.e.$  on  $\Omega \backslash D_{\epsilon
 }$  i.e.  $\displaystyle \Supp u\subset D_{\epsilon }.$ \ \par 
So we get\ \par 
\quad \quad \quad  $\forall \varphi \in {\mathcal{D}}_{n-p,n-q}(\Omega ),\ (-1)^{p+q-1}{\left\langle{\omega
 ,\ \varphi }\right\rangle}={\left\langle{u_{k},\ \bar \partial
 \varphi }\right\rangle}\rightarrow {\left\langle{u,\ \bar \partial
 \varphi }\right\rangle}\Rightarrow {\left\langle{u,\ \bar \partial
 \varphi }\right\rangle}=(-1)^{p+q-1}{\left\langle{\omega ,\
 \varphi }\right\rangle}$ \ \par 
hence  $\bar \partial u=\omega $  as distributions.\ \par 
\quad We proved that :\ \par 
\quad \quad \quad $\displaystyle \exists C>0,\ \forall \epsilon >0,\ \exists u_{\epsilon
 }\in L^{2}_{(p,q-1)}(\Omega ,c_{\varphi }e^{\varphi })::\bar
 \partial u_{\epsilon }=\omega ,\ {\left\Vert{u_{\epsilon }}\right\Vert}_{L^{2}(\Omega
 ,c_{\varphi }e^{\varphi })}\leq C{\left\Vert{\omega }\right\Vert}_{L^{2}(\Omega
 ,e^{\varphi })}^{2},$ \ \par 
and  $\displaystyle \Supp \ u_{\epsilon }\subset D_{\epsilon
 }.$  Because  $C$  is independent of  $\epsilon ,$  we can take
 again a subsequence to get finally that \ \par 
\quad \quad \quad $\displaystyle \exists u\in L^{2}_{(p,q-1)}(\Omega ,c_{\varphi
 }e^{\varphi })::\bar \partial u=\omega ,\ {\left\Vert{u}\right\Vert}_{L^{2}(\Omega
 ,c_{\varphi }e^{\varphi })}\leq C{\left\Vert{\omega }\right\Vert}_{L^{2}(\Omega
 ,e^{\varphi })}^{2},\ \Supp \ u\subset \bar D,$ \ \par 
and the theorem is proved.  $\displaystyle \hfill\blacksquare $ \ \par 
\quad We have a finer control of the support (see~\cite{AnGrauLr12},
 theorem 2.13).\ \par 

\begin{Thrm}
Let  $\displaystyle \Omega $  be a pseudo convex domain in  ${\mathbb{C}}^{n}$
 and   $\varphi $  be a strictly plurisubharmonic function of
 class  ${\mathcal{C}}^{2}$  in  $\displaystyle \Omega .$  Let
  $\displaystyle c_{\varphi }(z)$  be the smallest eigenvalue
 of  $\partial \bar \partial \varphi (z).$ \par 
Let  $\displaystyle \omega \in L^{2,c}_{p,q}(\Omega ,\ e^{\varphi
 })$  with  $\displaystyle \ \bar \partial \omega =0$  if  $\displaystyle
 q<n$  and   $\displaystyle \omega \perp {\mathcal{H}}_{n-p}(\Omega
 ,\ e^{-\varphi })$  if  $\displaystyle q=n,$  and with support
 in  $\Omega \backslash C$  where  $C$  is also a pseudo convex
 domain, then there is a  $\displaystyle u\in L^{2,c}_{p,n-1}(\Omega
 ,c_{\varphi }e^{\varphi }),\ \bar \partial u=\omega ,\ {\left\Vert{u}\right\Vert}\leq
 C{\left\Vert{\omega }\right\Vert}$  with support in  $\Omega
 \backslash \bar U,$  where  $U$  is any open set relatively
 compact in  $C,$  provided that  $q\geq 2.$  
\end{Thrm}
\quad \quad 	Proof.\ \par 
Let  $\omega $  be a  $(p,q)$  form with compact support in 
 $\Omega \backslash C$  then $\ \exists v\in L^{2,c}_{p,q-1}(\Omega
 ,c_{\varphi }e^{\varphi }),\ \bar \partial v=\omega ,$  with
 compact support in  $\Omega ,$  by theorem~\ref{cS1}.\ \par 
Because  $\omega $  has its support outside  $C$  we have  $\omega
 =0$  in  $C\ ;$  this means that  $\bar \partial v=0$  in  $C.$
  On the support of  $v$  we have, because  $\varphi $  is s.p.s.h.
 and of class  $\displaystyle {\mathcal{C}}^{2},$ \ \par 
\quad \quad \quad \begin{equation}  0<\delta \leq c_{\varphi }e^{\varphi }\leq
 M<\infty \label{cS5}\end{equation}\ \par 
hence  $\displaystyle v\in L^{2}_{p,q-1}(C).$  By H\"ormander
 we have, because  $C$  is Stein and bounded, that\ \par 
\quad \quad \quad $\displaystyle \exists h\in L^{2}_{p,q-2}(C)::\bar \partial h=v$
  in  $\displaystyle C.$ \ \par 
Let  $U$  be open and such that  $\displaystyle \bar U\subset
 C.$  Let  $\chi $  be a smooth function such that  $\chi =1$
  in  $U$  and  $\chi =0$  near  $\partial C\ ;$  then set\ \par 
\quad \quad \quad \quad \quad  $u:=v-\bar \partial (\chi h).$ \ \par 
We have that  $u=v-\chi \bar \partial h-\bar \partial \chi \wedge
 h=v-\chi v-\bar \partial \chi \wedge h$  hence  $u$  is in 
 $\displaystyle L^{2}_{p,q-1}(\Omega )\ ;$  moreover  $u=0$ 
 in  $\bar U$  because  $\chi =1$  in  $U$  hence  $\bar \partial
 \chi =0$  there. Finally  $\bar \partial u=\bar \partial v-\bar
 \partial ^{2}(\chi h)=\omega .$  Because of~(\ref{cS5})  $\displaystyle
 u\in L^{2}_{p,q-1}(\Omega )$  implies  $\displaystyle u\in L^{2}_{p,q-1}(\Omega
 ,\ c_{\varphi }e^{\varphi })$  and we are done.  $\hfill\blacksquare $ \ \par 
\quad \quad 	If  $\Omega $  and  $C$  are, for instance, pseudo-convex in
  ${\mathbb{C}}^{n}$  then  $\Omega \backslash C$  is no longer
 pseudo-convex in general, so this theorem improves actually
 the control of the support.\ \par 

\section{Approximation procedure.}

\begin{Thrm}
Let  $\displaystyle \omega \in L^{2}_{p,q}({\mathbb{C}}^{n},\
 e^{\varphi })$  with  $\displaystyle \bar \partial \omega =0$
  if  $\displaystyle q<n$  and   $\displaystyle \omega \perp
 {\mathcal{H}}_{n-p}(\Omega ,\ e^{-\varphi })$  if  $\displaystyle
 q=n,$  then there is  $\displaystyle u\in L^{2}_{p,q-1}({\mathbb{C}}^{n},c_{\varphi
 }e^{\varphi })$  such that  $\displaystyle \bar \partial u=\omega ,$  and\par 
\quad \quad \quad \quad \quad \quad \quad  $\displaystyle \ {\left\Vert{u}\right\Vert}_{L^{2}({\mathbb{C}}^{n},c_{\varphi
 }e^{\varphi })}\leq C{\left\Vert{\omega }\right\Vert}_{L^{2}({\mathbb{C}}^{n},e^{\varphi
 })}.$ 
\end{Thrm}
\quad This is theorem 1.2 of~\cite{AmarHeden13} where the assumption
 of the compact support for  $\omega $  is removed.\ \par 
\quad Proof.\ \par 
The  idea is to proceed by approximation. So by theorem~\ref{cS1}
 we have the following result:\ \par 
\quad \quad \begin{equation}  \forall \omega \in L^{2,c}_{p,q}({\mathbb{C}}^{n},\
 e^{\varphi }),\ \exists u\in L^{2,c}_{p,q-1}({\mathbb{C}}^{n},c_{\varphi
 }e^{\varphi }),\ \bar \partial u=\omega ,\ {\left\Vert{u}\right\Vert}\leq
 C{\left\Vert{\omega }\right\Vert}.\label{cS7}\end{equation}\ \par 
Now let  $\displaystyle \omega \in L^{2}_{p,q}({\mathbb{C}}^{n},\
 e^{\varphi })$  then there exists  $\displaystyle B_{k}=B(0,r_{k}),\
 r_{k}\geq r_{k-1}+1\rightarrow \infty ,$  a sequence of  balls
 in  $\displaystyle {\mathbb{C}}^{n}$   such that\ \par 
\quad \quad \quad \begin{equation}  \ \int_{B_{k}^{c}}{\left\vert{\omega }\right\vert
 ^{2}e^{\varphi }dm}\leq 1/k+1.\label{c1}\end{equation}\ \par 
Take  $\displaystyle \chi _{k}\in {\mathcal{C}}^{\infty }_{c}(B_{k})$
   such that  $\displaystyle \chi _{k}=1$   in  $\displaystyle
 B_{k-1}.$    Consider the form  $\displaystyle \omega _{k}:=\chi
 _{k}\omega $  ; we have:  $\displaystyle \ {\left\Vert{\omega
 -\omega _{k}}\right\Vert}\leq 1/k$  by~(\ref{c1}).\ \par 
And  $\displaystyle \bar \partial \omega _{k}=\bar \partial \chi
 _{k}\wedge \omega .$  Clearly we can choose  $\chi _{k}$  such
 that  $\displaystyle \ {\left\Vert{\bar \partial \chi _{k}}\right\Vert}_{\infty
 }\leq 1$  and this is the place where  $\displaystyle {\mathbb{C}}^{n}$
  is required.\ \par 
Hence\ \par 
\quad \quad \quad \quad $\displaystyle \ {\left\Vert{\bar \partial \omega _{k}}\right\Vert}^{2}=\int_{B_{k}\backslash
 B_{k-1}}{\left\vert{\bar \partial \chi _{k}}\right\vert ^{2}\left\vert{\omega
 }\right\vert ^{2}e^{\varphi }}\leq \int_{B_{k}\backslash B_{k-1}}{\left\vert{\omega
 }\right\vert ^{2}e^{\varphi }}\leq \frac{1}{k}.$ \ \par 
By~(\ref{cS7}) we  have that there exists  $\displaystyle u_{k}\in
 L^{2,c}_{p,q}({\mathbb{C}}^{n},c_{\varphi }e^{\varphi }),\ \bar
 \partial u_{k}=\bar \partial \omega _{k}$   with  $\displaystyle
 \ {\left\Vert{u_{k}}\right\Vert}\leq C/k.$ \ \par 
\quad Now set:  $\displaystyle \mu _{k}:=\omega _{k}-u_{k}$  ; we have
  $\displaystyle \bar \partial \mu _{k}=0$   and  $\displaystyle
 \ {\left\Vert{\omega -\mu _{k}}\right\Vert}\leq \frac{1+C}{k}$
  because  $\omega -\mu _{k}=\omega -\omega _{k}+u_{k}.$  Moreover
 we have that  $\displaystyle \mu _{k}$   is compactly supported.
 So there is  $\displaystyle v_{k}\in L^{2,c}_{p,q-1}({\mathbb{C}}^{n},c_{\varphi
 }e^{\varphi }),\ \bar \partial v_{k}=\mu _{k}$  and  $\displaystyle
 \ {\left\Vert{v_{k}}\right\Vert}\leq C{\left\Vert{\mu _{k}}\right\Vert}\leq
 C({\left\Vert{\omega }\right\Vert}+\frac{1+C}{k}).$ \ \par 
Because the norm of  $\displaystyle v_{k}$   in  $\displaystyle
 L^{2}_{p,q-1}({\mathbb{C}}^{n},c_{\varphi }e^{\varphi })$  
 is uniformly bounded, we have a weakly converging  sub sequence :\ \par 
\quad \quad \quad $\displaystyle \exists v\in L^{2}_{p,q-1}({\mathbb{C}}^{n},c_{\varphi
 }e^{\varphi })::\forall f\in L^{2}_{n-p,n-q+1}({\mathbb{C}}^{n},c_{\varphi
 }e^{\varphi }),\ \lim \ {\left\langle{v_{k},f}\right\rangle}={\left\langle{v,f}\right\rangle}.$
 \ \par 
Now choose  $\displaystyle f:=\bar \partial \varphi ,$  with
  $\displaystyle \varphi \in {\mathcal{D}}_{n-p,n-q}({\mathbb{C}}^{n})$
  then\ \par 
\quad \quad \quad $\displaystyle \ {\left\langle{v_{k},\bar \partial \varphi }\right\rangle}=(-1)^{p+q}{\left\langle{\bar
 \partial v_{k},\varphi }\right\rangle}=(-1)^{p+q}{\left\langle{\mu
 _{k},\varphi }\right\rangle}\rightarrow {\left\langle{v,\bar
 \partial \varphi }\right\rangle}=(-1)^{p+q}{\left\langle{\bar
 \partial v,\varphi }\right\rangle}.$ \ \par 
But  $\displaystyle \mu _{k}=\omega +e_{k},$  with  $\displaystyle
 \ {\left\Vert{e_{k}}\right\Vert}\leq \frac{1+C}{k}$  hence\ \par 
\quad \quad \quad $\displaystyle \ {\left\langle{\mu _{k},\varphi }\right\rangle}={\left\langle{\omega
 ,\varphi }\right\rangle}+{\left\langle{e_{k},\varphi }\right\rangle}\rightarrow
 {\left\langle{\omega ,\varphi }\right\rangle}$ \ \par 
so\ \par 
\quad \quad \quad $\displaystyle \ {\left\langle{\bar \partial v,\varphi }\right\rangle}={\left\langle{\omega
 ,\varphi }\right\rangle}$ \ \par 
which means that  $\displaystyle \bar \partial v=\omega $  in
 the distributions sense.  $\displaystyle \hfill\blacksquare $ \ \par 
\ \par 

\bibliographystyle{/usr/local/texlive/2013/texmf-dist/bibtex/bst/base/plain}

\end{document}